# Computational Approaches for Solving Two-Echelon Vehicle and UAV Routing Problems for Post-Disaster Humanitarian Operations

**Tasnim Ibn Faiz[1], Chrysafis Vogiatzis[2], and Md. Noor-E-Alam[1*]**

[1]Department of Mechanical and Industrial Engineering,
Northeastern University
360 Huntington Ave, Boston, MA 02115, USA

[2]Department of Industrial and Enterprise Systems Engineering
University of Illinois at Urbana-Champaign
104 S. Mathews Ave, Urbana, IL 61801, USA

*Corresponding author email: mnalam@neu.edu

## ABSTRACT

Humanitarian logistics service providers have two major responsibilities immediately after a disaster: locating trapped people and routing aid to them. These difficult operations are further hindered by failures in the transportation and telecommunications networks, which are often rendered unusable by the disaster at hand. In this work, we propose a two-echelon vehicle routing framework for performing these operations using aerial uncrewed autonomous vehicles (UAVs or drones) to address the issues associated with these failures. In our proposed framework, we assume that ground vehicles cannot reach the trapped population directly, but they can only transport drones from a depot to some intermediate locations. The drones launched from these locations serve to both identify demands for medical and other aids (e.g., epi-pens, medical supplies, dry food, water) and make deliveries to satisfy them. Specifically, we present a decision framework, in which the resulting optimization problem is formulated as a two-echelon vehicle routing problem with trucks as the first echelon vehicles and for the second echelon vehicles we consider two types of drones. Hotspot drones have the capability of providing cell phone and internet reception, and hence are used to capture demands. Delivery drones are subsequently employed to satisfy the observed demand. To handle demand uncertainty, we decompose the decision problem into two stages: providing telecommunications capabilities in the first stage thereby capturing demand precisely and satisfying the resulting demands in the second stage. To solve the resulting models, we propose efficient computational approaches by designing a decomposition algorithm with column generation (CG)-based heuristics to identify optimal drone routes. To showcase the applicability of our proposed framework, we present results from the numerical experiments on datasets created to simulate the demand for medical aid in Puerto Rico after Hurricane Maria.

## 1. Introduction

In the aftermath of a disaster, locating and sending assistance to people that are trapped inside the disaster zone is of utmost importance. Trapped populations are unable to leave their positions due to infrastructure failures (e.g., flooded streets after hurricanes, destroyed bridges after earthquakes), and are often unable to communicate their exact locations and needs, due to widespread communication and power systems failures. Even if humanitarian aid agencies become aware of the trapped population pressing needs for medical supplies, potable water, and dry food, they cannot send help via ground vehicles due to the potentially unusable transportation network following the disaster.

A real-life motivating example is the disastrous impact of Hurricane Maria in Puerto Rico in 2017. After the hurricane, much of the island of Puerto Rico was left without functional communication infrastructure. Many people were trapped in their communities with limited ability to move and without means for communicating their need for food and medical supplies. Furthermore, road transportation networks had failed and reaching the people in need using ground vehicles was impossible. In cases like these, uncrewed aerial vehicles (UAVs), assisted and transported by ground vehicles, can be used to (i) provide telecommunication capabilities to the affected areas, and (ii) deliver emergency medical and food supplies to the people in need.

In our work, we develop mathematical models for scheduling and routing ground and aerial vehicles with the double goals of efficiently providing telecommunications capabilities and delivering emergency supplies. The problem can be viewed as a variant of the two-echelon capacitated vehicle routing problem, in which ground vehicles (trucks, or first echelon vehicles) are tasked with carrying a finite number of uncrewed aerial vehicles (UAVs, drones, or second echelon vehicles) equipped with limited medical and other emergency supplies. Due to infrastructure failures and to ensure the wellbeing and safety of first responders and emergency personnel, ground vehicles can only reach so-called satellite locations outside the disaster zones, with drones being deployed to carry out the delivery operations within the zones themselves.

More specifically, in the proposed framework (referred to hereafter as **2EVRP-HD-UAV**), we consider two types of uncrewed aerial vehicles: (i) hotspot drones that are responsible for first restoring communication signals in an area and can then receive the demand requests from people in the area, and (ii) delivery drones that drop off medical supplies at the affected communities (demand locations). The demand information becomes available only after the communication signal is provided at all the demand locations by the hotspot drones. The demand information then dictates the delivery route decisions. We design a two-stage optimization approach to handle the demand materialization after the hotspot drone operations. In this framework, the first echelon ground vehicle routing and scheduling decisions are also made in tandem with the UAV routing decisions to minimize the time to reach all affected communities with aid deliveries and the cost of failure to satisfy demand. This setup naturally results in idle times for the ground vehicles at the satellite stations due to hotspot drone operations but results in precise demand information.

This paper is outlined as follows. First, we provide the contributions of the work to the body of literature. We then provide an accounting of the relevant literature in Section 2. As our work is in the intersection of multiple vehicle routing frameworks, we present a general overview of two-echelon vehicle routing problems, related solution approaches, and their applications involving uncrewed autonomous vehicles (aerial or otherwise). In Section 3, we present our framework with the mathematical formulation to model it, and the solution approach to solve it. We finish this work with our computational results in Section 4. There, we provide a case study on a simulated dataset for relief material demand at the zip code areas in Puerto Rico immediately after Hurricane Maria. This is a prime example for our setup as multiple counties in Puerto Rico were left without power and communications for many days after the hurricane made landfall. We also show the validity of our model on a synthetic dataset (presented in Appendix A2). Finally, we provide our concluding remarks in Section 5.

### *1.1. Contributions*

Considering the state of the art (see Section 2 and the detailed literature review), we see that exact approaches have not been particularly successful in real-world two-echelon vehicle routing problems (2E-VRP) with and without uncrewed autonomous vehicles (UAVs). The scale and complexity of the aid routing problems usually encountered in disaster management makes the use of exact approaches difficult. To overcome the limitations of these existing approaches, our study contributes to the development of models and algorithms for two-echelon vehicle routing problems using drones in a post-disaster environment where demand quantities at the affected communities are uncertain. We specifically propose a **two-stage decision approach to handle demand uncertainty**. More specifically, our contributions can be summarized as follows.

First, we investigate a novel 2E-VRP with trucks carrying UAVs as the first echelon vehicles and the UAVs as the second echelon vehicles for performing humanitarian logistics operations. Our fundamental underlying assumption and driver for this research is that at the time of truck and UAV route determination, the actual aid demand at the affected communities is not known. For this novel problem, termed 2EVRP-HD-UAV (two echelon vehicle routing problem with hotspot and delivery UAVs), we present a mixed integer linear program (MILP). Unsurprisingly, the MILP is computationally expensive to solve even for smaller size instances (e.g., ones with 12 hotspot locations, 20 demand locations).

Second, to that end, we present a set covering formulation, which is suitable for column generation. We then propose a computational framework for heuristic solution approaches that are based on column generation (CG). Our decomposition approach leveraging the problem structure enables us to set up very small pricing subproblems for each satellite-drone pair to gain computational efficiency. The proposed approach can solve problems of up to 120 demand locations within reasonable computational time, as shown in our experiments.

Finally, to handle demand uncertainty, we design a two-stage decision framework. In the first stage, we identify the optimal routes for trucks, hotspot, and delivery drones to provide coverage to all demand location, while the second stage decisions involve identifying the best delivery routes based on the demand information. The assumption is that the demand information becomes available only after the (first stage) routing decisions are made. We formulate the second stage model fixing the truck and hotspot drone routing decisions found at the termination of the first stage. Following a CG-based heuristic (a matheuristic, or math-based heuristic), we generate new delivery drone routes with the demand information if some of the routes generated in the first stage become infeasible. In our experiments, we show that our CG-based computational approaches perform better for the proposed two-stage decision framework compared to exact methods for solving the MILP model.

## 2. Literature review

In this section, we provide a general overview of the relevant literature. First, we briefly mention recent work in vehicle routing problems in the context of humanitarian logistics and disaster management. We then turn our focus to two-echelon vehicle routing problems and provide a survey of the state-of-the-art exact and heuristic solution approaches available. Moreover, we discuss studies that involve routing of UAVs (drones) along with ground vehicles for optimizing deliveries. Finally, we present work on robust vehicle routing problems with uncertain data (as in demand uncertainty).

We begin with a brief discussion of recent work in vehicle routing problems in humanitarian aid and relief operations. Interested readers are referred to the review works in [1-3], that cover literature which largely inspired our problem definition and our assumptions. Vehicle routing problems (VRP) have been studied extensively in this context of natural and human-made disasters [4-8]. The underlying assumption of these studies is that of a fully functional post-disaster transportation network, which implies that the methods and ideas are not always adoptable in cases where the supply chain infrastructure itself gets affected by the disaster. In the absence of supply chain infrastructure, the adoption of new technologies and approaches are necessary for fulfilling aid demand, which is indeed the premise of our study.

Two-echelon vehicle routing problems (2E-VRP) are extensions of the classical VRP, where the first echelon connects the depot to the intermediate satellite locations and the second echelon connects the satellite locations to the customers. Typically, the objective is to minimize the cost associated with the operations in both echelons. An excellent survey outlining the state-of-the-art advancements in the area can be found in [9]. Exact solution methods for the capacitated version of 2E-VRP (commonly referred to as 2E-CVRP) are primarily based on branch-and-cut algorithms, originally proposed by Feliu et al. [10] and improved in subsequent studies [11-13]. Baldacci et al. [14], however, employ a hybrid exact algorithm to decompose the problem into a set of capacitated vehicle routing problems. Combining the strengths of branch-and-cut and branch-and-price, Santos et al. [15, 16] use a branch-and-cut-and-price algorithm, strengthened with several classes of valid inequalities. More recently, strong valid inequalities within a branch-and-cut scheme have received significant attention [17, 18].

Due to the computational complexity associated with solving problems of the 2E-CVRP family using exact methods, heuristic approaches have gained traction. Clustering-based heuristics [19] and multi-start heuristics [20] showcase the impact of satellite locations on the cost improvement opportunities in 2E-CVRP. An adaptive large neighbourhood search (ALNS) [21] algorithm, a hybrid metaheuristic based on large neighbourhood search (LNS) [22], and a hybrid multi-population genetic algorithm [23] have also been used to solve 2E-CVRP problems. Finally, a greedy randomized adaptive search procedure (GRASP) combined with a route-first cluster-second and a variable neighbourhood descent algorithm has also been proposed in [24].

A specific aspect of 2E-VRP involves cooperation or synchronization between trucks (first echelon) and UAVs/drones (second echelon). In this setup, UAVs are launched from and land back to a truck to recharge batteries, to be reloaded with delivery items, or to be transported to the next launching location. Review works on VRPs with synchronization and cooperative truck-drone decision problems can be found in [25, 26]. In the next paragraphs, we discuss some of the recent works in this area.

A two-echelon cooperative truck and drone routing problem is solved by a branch-and-cut framework with a separation procedure [27], whereas heuristic approaches are used in [28]. Wang et al. [29], on the other hand, derive several worst-case results for multiple trucks and drones problems. The study has later been extended in [30, 31], where drones are allowed to provide close-enough coverage to the demand points. Carlsson and Song [32] follow a continuous approximation approach for solving this family of problems.

Recently, last-mile deliveries via truck and drone collaboration have resulted in a new variant of the traveling salesperson problem (TSP), termed TSP with drone (TSP-D). Murray and Chu study the simultaneous routing of a truck and a drone [33] and later Agatz et al. propose several route-first cluster-second heuristics [34] to solve these problems. A decomposition-based iterative approach is also used to solve TSP-D problems [35]. Ha et al. [36] propose heuristic approaches to solve the TSP-D problem, whereas a combination of *k*-means clustering and TSP is proposed by Chang and Lee [37]. A case with multiple trucks and drones is studied in [38], where a multiple traveling salesperson problem with drones (mTSPD) is solved using an adaptive insertion heuristic. A vehicle-assisted multi-drone routing and scheduling problem is studied in [39], with an algorithm based on the adaptive memory procedure is proposed to solve it. A route optimization method using genetic algorithm for a moving truck and a drone launching problem is presented by Savuran and Karakaya [40, 41].

The existing literature on routing UAVs has been mostly focused on deterministic variants. For example, in a recent work, Ghelichi et al. [42] proposed a model for locating charging stations and scheduling a fleet of drones for delivery of medical items. Stochastic extensions are an exception: for example, a robust optimization approach is proposed by Evers et al. [43] to plan surveillance missions using UAVs subject to uncertain weather conditions. A similar approach is used by Kim et al. [44] to address uncertainty in battery capacity due to air temperature. In contrast, literature on 2E-VRP with stochastic parameters is limited. One of these studies adopts a genetic algorithm based method [45]. The other research work employs a simulation-based optimization approach [46], in which a Monte Carlo sampling method is used to handle stochastic demand.

## 3. The 2EVRP-HD-UAV framework

In this section, we present the two-stage decision framework for 2E-VRP with hotspot and delivery UAVs followed by a CG-based two-stage solution heuristic. Following a disaster, in the absence of adequate transportation network infrastructure, affected communities are not directly reachable by trucks; in addition, in the absence of communications infrastructure, the demand levels at the communities are not known precisely. We discuss the necessary assumptions used to formulate the models for the 2EVRP-HD-UAV framework below:

1. All the trucks and the hotspot and delivery drones are homogeneous (that is, all trucks are identical, and all drones are identical).
2. The weather conditions and vehicle characteristics are deterministic. Vehicle travel speeds and other details are estimated considering the specific post-disaster operating conditions.

3. The trucks that carry hotspot and delivery drones are dispatched from a central depot and can visit satellite stations in sequence. Satellite stations are selected from a suitable set of candidate locations in major highways that are still usable after a disaster. A satellite location can be visited by at most one truck.
4. Hotspot and delivery drones are launched from a truck that is stopped at a satellite and are routed to one or more hotspot locations and one or more affected communities, respectively before they return to that same satellite.
5. Each affected community has at least one satellite location in its vicinity so that a hotspot drone can cover it. We further assume that each community has at least one functioning device to send information to and from.
6. Once a hotspot drone reaches a hotspot location, it hovers until all the affected communities that are within a certain distance receive telecommunication signal: this is how the demand information of these communities is obtained.
7. The velocity of hotspot drones is estimated considering the time for hovering to setup telecommunication signal and capture demand information.
8. The times for delivery drones are estimated considering the times required for loading, stopping, and delivering aid at disaster affected communities.
9. Delivery drones can be launched from a satellite only after the hotspot drones launched from the same satellite return after their tours.
10. A delivery drone can be dispatched from a satellite to deliver aid only if the community has been provided with coverage from a hotspot drone also launched from that same satellite.
11. The truck and hotspot drone routing decisions are made in a centralized manner to provide telecommunication network and delivery coverage to all communities. The demand information gathered by the hotspot drones impacts the delivery drone operations at a local (satellite) level. However, it does not impact the truck routing decisions.

With these assumptions at hand, we provide a more formal problem definition in the following section.

### *3.1. Problem definition*

To formulate the problem, we consider a graph $G(V, E)$. Node set $V = \{0\} \cup S \cup B \cup C$ consists of the depot (0), the set of candidate satellite locations ($S$), the set of candidate hotspot locations ($B$), and a set of affected communities ($C$). The edge set $E$ consists of affected and unaffected edges, i.e., $E = E^a \cup E^o$. Drones can use any edge in the network, unconstrained from the infrastructure condition (e.g., flooded roads, destroyed bridges); on the other hand, ground vehicles can only use edges in $E^o$. We assume that a fleet of identical trucks are available that can transport a fixed number of hotspot and delivery drones. For all edges that are operational ($\forall(i,j) \in E^o, i,j \in \{0\} \cup S$), we denote traversal times for trucks as $\tau_{ij}^t, \forall(i,j) \in E^o$. The affected edge set $E^a$, on the other hand, are comprised of two subsets of edges, $E^h$ and $E^d$ which are available for traversing by hotspot and delivery drones, respectively. The edge set $E_s^h$ is the collection of edges that can be traversed by a hotspot drone launched from a truck that stops and stays at satellite location $s \in S$, i.e., $\forall(i,j) \in E_s^h, i,j \in \text{s} \cup B_s$, where $B_s$ is the set of hotspot nodes located within reachable distance from satellite s; if a hotspot drone is launched from $s$, the total roundtrip time from $s$ to $b$ is less than the hotspot drone flying time range. The edge set $E^h$ for hotspot drone traversing is composed of subset of edges $E_s^h, \forall s \in S$, i.e., $E^h = \bigcup_{s \in S} E_s^h$. Edge traversing times for hotspot drones are given by $\tau_{ij}^h, \forall(i,j) \in E^h : i,j \in S \cup B$. Similarly, the delivery edge subset $E^d$ consists of smaller subsets $E_s^d, \forall s \in S$, i.e., $E^d = \bigcup_{s \in S} E_s^d$. Arcs in $E_s^d$ are composed of satellite node $s$ and the nearby community locations $C_s$ that are reachable by a delivery drone launched from $s$. Edge traversing times for delivery drone are $\tau_{ij}^d, \forall(i,j) \in E^d : i,j \in S \cup C$. The hotspot drones are dispatched and routed so that they visit hotspot locations, and a drone visiting a hotspot node $b \in B$ can provide telecommunication signals to all the affected communities that are within a specified distance from $b$. Once a hotspot drone finishes its route, it returns to its launching satellite location. Only then the delivery drones can be dispatched from the same satellite location to satisfy demands at the communities that have been covered by hotspot drones from the same satellite. In Figure 1, we show the pictorial representation of the decision problem. For better visualization, we only show a few vehicle routes. Red arcs present the truck routes, whereas hotspot and delivery drone routes are shown in green and yellow arcs respectively.

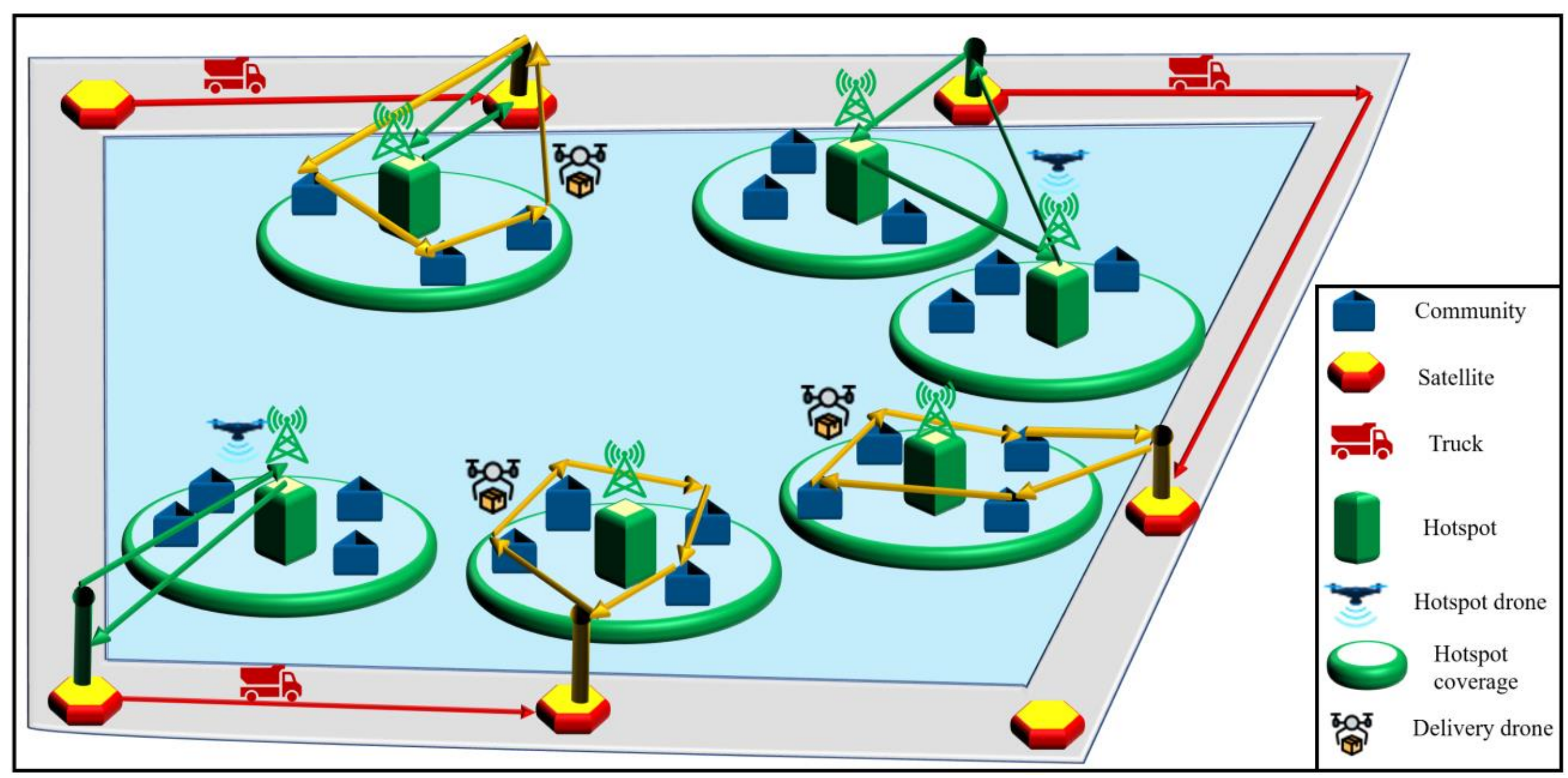

Figure 1: A pictorial representation of the decision problem*

The first stage decision problem of the 2EVRP-HD-UAV consists of routing each of the available trucks to one or more satellite locations, dispatching and routing of hotspot drones to provide telecommunication signal to all communities, and deploying and routing delivery drones to the affected communities. The objective is to minimize the cost of delay in reaching all communities, along with a (penalty) cost of failing to reach a community by delivery drones.

To solve the resulting 2E-VRP with two classes of second echelon vehicles, we first develop a mixed integer linear programming model and find that the problem is too computationally expensive to solve even smaller size instances. For completeness, we present this model in the Appendix. That said, in the next subsections we present a set covering problem reformulation that enables us to develop a column generation-based decomposition and solution approach. We begin by presenting the notations used in the model.

* Drone icons obtained from Pixel perfect and www.flaticon.com.

### *3.2. First stage set covering problem formulation*

All the necessary notation and the associated parameters used in the remainder of the section are presented in the following tables.

*Sets:*

| | |
|---|---|
| $S$ | Set of all satellite locations |
| $B$ | Set of all locations permitted for hotspot drone hovering; $B = \bigcup_{s \in S} B_s$ |
| $C$ | Set of all community locations; $C = \bigcup_{s \in S} C_s$ |
| $H_s$ | Set of hotspot drones in a truck at satellite $s \in S$ |
| $D_s$ | Set of delivery drones in a truck at satellite $s \in S$ |
| $E^o$ | Set of arcs for trucks; $(i,j) \in E^o: i,j \in S \cup \{0\}$ |
| $E_s^h$ | Set of arcs for hotspot drones starting from satellite $s \in S$; $(i,j) \in E_s^h: i,j \in s \cup B$ |
| $E_s^d$ | Set of arcs for delivery drones starting from satellite $s \in S$; $(i,j) \in E_s^d: i,j \in s \cup C$ |
| $U_s$ | Restricted set of hotspot drone routes starting from satellite $s \in S$ |
| $V_s$ | Restricted set of delivery drone routes starting from satellite $s \in S$ |

*Parameters:*

| | |
|---|---|
| $m^t$ | Number of trucks available |
| $m^h$ | Number of hotspot drones carried in a truck |
| $m^d$ | Number of delivery drones carried in a truck |
| $F_c^T$ | Cost of unit time delay in reaching community $c \in C$ by a delivery drone |
| $F_c^R$ | Cost of failing to reach community $c \in C$ by a delivery drone |
| $\tau_{ij}^t$ | Time to traverse arc $(i,j) \in E^o$ by a truck |
| $\tau_{ij}^h$ | Time to traverse arc $(i,j) \in E_s^h$ by a hotspot drone |
| $\tau_{ij}^d$ | Time to traverse arc $(i,j) \in E_s^d$ by a delivery drone |
| $W^h$ | Flying range of a hotspot drone |
| $W^d$ | Flying range of a delivery drone |
| $\Pi_{sk}^h$ | Route time length for route $k \in U_s$ of a hotspot drone from satellite $s \in S$ |
| $\Pi_{sp}^d$ | Route time length for route $p \in V_s$ of a delivery drone from satellite $s \in S$ |
| $t_{sc}^k$ | Time to reach community $c \in C$ by delivery drone $k \in D_s$ from satellite $s \in S$ |
| $a_j^k$ | =1, if hotspot location $j \in B$ is present in $k \in U_s$ |
| $b_c^k$ | =1, if community $c \in C$ is present in $k \in V_s$ |
| $g_{jc}$ | =1, if demand at community $c \in C$ can be confirmed by a hotspot drone hovering at $j \in B$, 0 otherwise |

*Decision variables:*

| | |
|---|---|
| $x_{ij}$ | =1, if a truck traverses arc $(i,j) \in E^o$; 0 otherwise |
| $y_p^{sk}$ | =1, if hotspot drone $k \in H_s$ travels route $p \in U_s$ starting from $s \in S$; 0 otherwise |
| $z_p^{sl}$ | =1, if delivery drone $l \in D_s$ travels route $p \in V_s$ starting from $s \in S$; 0 otherwise |
| $J_c$ | =1, if community $c \in C$ cannot be reached by a delivery drone |
| $\mathrm{T}_{sc}^d$ | Time to reach community $c \in C$ by a delivery drone launched from satellite $s \in S$ |
| $\Delta_s^h$ | Length of time spent by a truck at satellite $s \in S$ due to hotspot drone flights |
| $\Delta_s^d$ | Length of time spent by a truck at satellite $s \in S$ due to delivery drone flights |
| $T_s$ | Arrival time of a truck at satellite $s \in S$ |

The set covering formulation of the first stage enables us to develop a column generation-based decomposition heuristic, where the problem is decomposed into a main problem and its subproblems. The restricted main problem of the first stage (D-RMP-1) can now be presented in (1)—(15). We call this problem restricted because we only consider restricted sets of routes for the hotspot and the delivery drones.

$$Minimize \sum_{s \in S} \sum_{c \in C} F_c^T \mathrm{T}_{sc}^d + \sum_{c \in C} F_c^R J_c \tag{1}$$

Subject to: **Duals**

$$\sum_{(i,j) \in E^o : i \notin S} x_{ij} \leq m^t \tag{2}$$

$$\sum_{(i,j) \in E^o} x_{ij} \leq 1 \qquad \forall j \in S, \tag{3}$$

$$\sum_{(i,j) \in E^o} x_{ij} - \sum_{(j,k) \in E^o} x_{ij} = 0 \qquad \forall j \in S, \tag{4}$$

$$\begin{aligned} T_j &\geq T_i + \Delta_i^h + \Delta_i^d + \tau_{ij}^t - M\left(1 - x_{ij}\right) && \forall (i,j) \in E^o \\ T_j &\leq T_i + \Delta_i^h + \Delta_i^d + \tau_{ij}^t + M\left(1 - x_{ij}\right) && \forall (i,j) \in E^o \end{aligned} \tag{5}$$

$$\sum_{k\in H_s}\sum_{l\in U_s} y_l^{sk} - m^h \sum_{(i,s)\in E^o} x_{is} \le 0 \qquad \forall s\in S, \qquad \gamma_s \tag{6}$$

$$\sum_{k\in D_s}\sum_{l\in V_s} z_l^{sk} - m^d \sum_{(i,s)\in E^o} x_{is} \le 0 \qquad \forall s\in S, \qquad \theta_s \tag{7}$$

$$\sum_{l\in U_s} y_l^{sk} - \sum_{(i,s)\in E^o} x_{is} \le 0 \qquad \forall s\in S, k\in H_s \qquad \varphi_{sk} \tag{8}$$

$$\sum_{l\in V_s} z_l^{sk} - \sum_{(i,s)\in E^o} x_{is} \le 0 \qquad \forall s\in S, k\in D_s \qquad \phi_{sk} \tag{9}$$

$$\sum_{s\in S}\sum_{k\in D_s}\sum_{l\in V_s} z_l^{sk}\, b_{lc} + J_c = 1 \qquad \forall c\in C \qquad \zeta_c \tag{10}$$

$$\sum_{l\in U_s} y_l^{sk}\, \Pi_{sl}^h - \Delta_s^h \le 0 \qquad \forall s\in S, k\in H_s \qquad \sigma_{sk} \tag{11}$$

$$\sum_{l\in V_s} z_l^{sk}\, \Pi_{sl}^d - \Delta_s^d \le 0 \qquad \forall s\in S, k\in D_s \qquad \rho_{sk} \tag{12}$$

$$T_s + \Delta_s^h + \sum_{k\in D_s}\sum_{l\in V_s} z_l^{sk} t_{sc}^l - M\left(1 - \sum_{k\in D_s}\sum_{l\in V_s} z_l^{sk} b_c^l\right) \le \mathrm{T}_{sc}^d, \quad \forall s\in S, c\in C, \quad \alpha_{sc} \tag{13}$$

$$\sum_{l\in U_s}\sum_{j\in H_s}\sum_{p\in B:g_{pc}=1} a_p^l\, y_l^{sj} - \sum_{k\in D_s}\sum_{i\in V_s} z_i^{sk}\, b_c^i \ge 0 \qquad \forall s\in S, c\in C, \qquad \pi_{sc} \tag{14}$$

$$x_{ij}, y_m^{sk}, z_p^{sl}, J_c \in \{0,1\}; \quad \Delta_s^h, \Delta_s^d, T_s, \mathrm{T}_{sc}^d \ge 0$$
$$\forall s\in S, (i,j)\in E^o, m\in U_s, p\in V_s, k\in H_s, l\in D_s, c\in C \tag{15}$$

The objective, shown in (1), is to minimize a comprehensive cost function, which includes both cost and time components. More specifically, the objective function minimizes: (i) the total time to reach all community, and (ii) the penalty for failure to reach communities by delivery drones. We omit any fixed cost associated with truck or drone dispatch since we assume that a fixed set of trucks is available, each of them equipped with a fixed number of hotspot and delivery drones.

Constraints (2) limit the maximum number of trucks that can be used, while the number of trucks allowed to visit a satellite location is restricted to at most one by constraints (3). Flow preservation for trucks at satellite nodes is ensured by constraints (4). Constraints (5) dictate the time it takes a truck to reach each satellite node. In constraints (6) and (8), we ensure that hotspot drones are only dispatched from a satellite node if a truck has visited that node; furthermore, the number of hotspot drones dispatched is limited by the number of drones carried by that truck. A similar set of constraints, given in (7) and (9), applies to delivery drones. Constraints (10) ensure that every community is either visited by a delivery drone or the variable $J_c$ takes a non-zero value. According to constraints (11)—(12), the lengths of time that a truck has to spend at a satellite location due to flights of hotspot and delivery drones are bounded by the maximum flight times of hotspot and delivery routes originating from that satellite, respectively. Constraint set (13) defines the required times to reach every community by a delivery drone from each satellite. Constraints (14) are coverage constraints that allow for a community to be visited by a delivery drone dispatched from a satellite only if that location has been covered by a hotspot drone launched from the same satellite. Finally, variable restrictions are presented in (15).

The D-RMP-1 formulation enables us to consider only a subset of routes to solve (1)—(15) instead of enumerating every possible route for the hotspot and delivery drones. The generation of routes and determination of route parameters for hotspot and delivery drones are relegated to the pricing problems.

Since we have two different drone route sets in our problem: one for hotspot drones and the other for delivery drones, we initially attempted to develop a framework where a pricing problem generates routes for both types of drones simultaneously. But we found the resulting pricing problem is computationally very expensive. To formulate smaller and easier to solve pricing problems, we designed two decomposition approaches. In the first approach, we follow a satellite-level decomposition, where one pricing problem is formulated for each satellite. We present the pricing problems resulting from this decomposition scheme in the Appendix. In the second decomposition approach, we design satellite-drone-level decomposition approach, where one pricing problem is formulated for each satellite and each drone carried on a truck visiting the satellite. We found the second decomposition approach to be computationally more efficient. The pricing problem formulation resulting from the second decomposition approach for generating delivery drone routes is presented in Section 4.3. We also formulate a similar pricing problem for generating hotspot drone routes, which we present in the Appendix. We present the delivery drone route generation pricing problem below starting with the brief explanations of the notations.

Before showing the pricing problem, we present some more notation.

*Newly defined decision variables for the drone route generation pricing subproblem:*

$v_{ij}^{sk}$ =1, if arc $(i,j) \in \bigcup_{s \in S} E_s^d$ is traversed by delivery drone $k \in D_s$; 0 otherwise

$\mathrm{t}_j^{sk}$ Time to reach node $j \in C_s \cup s$ by delivery drone $k \in D_s$

### *3.3. Pricing subproblems for generating delivery drone routes (D-PSP-DD-1)*

The pricing subproblems for generating route for each delivery drone $k \in D_s$ from satellite $s \in S$, are given in (16)—(22).

$$\begin{aligned} \mathit{Minimize} \quad & (RC)_{sk}^d \\ & = -\theta_s - \phi_{sk} - \sum_{j \in C_s} \sum_{(i,j) \in E_s^d} v_{ij}^{sk} \left( \zeta_j + \pi_{sj} + \rho_{sk} \tau_{ij}^d \right) \\ & - \sum_{j \in C_s} \alpha_{sj} \left( \mathrm{t}_j^{sk} + M \sum_{(i,j) \in E_s^d} v_{ij}^{sk} \right) \end{aligned} \tag{16}$$

Subject to:

$$\sum_{(i,j) \in E_s^d} v_{ij}^{sk} \, \tau_{ij}^d \leq W^d \tag{17}$$

$$\sum_{(s,j)\in E_d^s} v_{sj}^{sk} \leq 1 \tag{18}$$

$$\sum_{(i,j)\in E_s^d} v_{ij}^{sk} \leq 1 \qquad \forall j \in C_s, \tag{19}$$

$$\sum_{(i,j)\in E_s^d} v_{ij}^{sk} - \sum_{(j,i)\in E_s^d} v_{ij}^{sk} = 0 \qquad \forall j \in C_s, \tag{20}$$

$$\begin{aligned} \mathrm{t}_j^{sk} &\geq \mathrm{t}_i^{sk} + \tau_{ij}^d - W^d\left(1 - v_{ij}^{sk}\right) && \forall (i,j) \in E_s^d \\ \mathrm{t}_j^{sk} &\leq \mathrm{t}_i^{sk} + \tau_{ij}^d + W^d\left(1 - v_{ij}^{sk}\right) && \forall (i,j) \in E_s^d, \end{aligned} \tag{21}$$

$$v_{ij}^{sk} \in \{0,1\}, \forall s \in \mathrm{S}, (i,j) \in E_s^d, k \in D_s, \quad \mathrm{t}_j^{sk} \geq 0, \ \forall j \in C_s, \tag{22}$$

The objective function in (16) finds the best reduced cost $(RC)_{sk}^d$ of delivery route for drone $k \in D_s$ launched from satellite $s \in S$. If the obtained reduced cost is negative, then that route is added to the restricted set $V_s$. Constraints (17) limits the total travel distance of a drone. Each drone can be used only along a single path starting from the satellite under consideration according to constraints (18). Each community location can be visited at most once per constraints (19); constraints (20) provide the balancing of inbound and outbound arcs coincident to a community location. Constraints (21) indicate the time to visit two consecutive nodes on delivery drone route. Finally, constraints (22) provide non-negativity bounds. As a reminder, $\theta_s$, $\phi_{sk}$, $\zeta_c$, $\rho_{sk}$, $\alpha_{sc}$, $\pi_{sc}$; $\forall s \in S, k \in D_s, c \in C_s$ are the dual variables associated with constraints (7), (9), (10), (12), (13), and (14).

It is necessary to mention here that since all the delivery drones are assumed identical, the best route for the first delivery drone aboard a truck stopped at satellite $s$ ($k \in D_s: k = 1$) would be considered as the best route for each of the other delivery drones ($k \in D_s: k = 2., \dots, m^d$) on that truck. The same circumstances are true for hotspot drones. To avoid generating the same route for all the drones from satellite $s$, we update the list of nodes to be visited ($B_s$ for hotspot drones and $C_s$ for delivery drones) for each pricing problem at every iteration. This update is done by dropping the nodes that are included in the already generated routes from the same iteration. For example, for the first delivery drone from the truck stopped at satellite $s$, it will consider all the community nodes that can be reached from the satellite ($C_s$) for route generation, whereas for the second drone, that is $k \in D_s: k = 2$, the set of nodes to be visited will be updated to $C_s \coloneqq C_s - \{j \in C_s: \sum_{(i,j)\in E_s^d} v_{ij}^{s1} = 1\}$ by removing the nodes that are in the first drone's route. This results in faster solution of pricing problems and distinct routes at the expense of (slightly) higher route costs.

The first stage model provides us with the routes and schedules for trucks, hotspot drones and delivery drones, which are then passed on to the second stage problem. The decision problem in this stage is to select optimal drone routes from the candidate route set to satisfy the demand at the affected communities while abiding by the number of drone and load carrying capacity constraints. If the existing route set fails to satisfy the demands, new delivery routes are generated while the truck and hotspot drone operations are assumed to be fixed and irreversible.

### *3.4. Second stage restricted main problem*

The solution from D-RMP-1 is used as input to the second stage model (D-RMP-2). The demands at the affected communities are known precisely in the second stage. With these, the second stage problem determines the delivery drone routes to satisfy the demand at the affected communities. The parameters and variables in this stage are defined below.

*Newly defined set and parameters:*

| | |
|---|---|
| $V_s'$ | Restricted set of delivery drone routes |
| $F_c^D$ | Cost of failing to satisfy each unit demand at community $c \in C$ |
| $F^E$ | Cost of exceeding carrying capacity of a delivery drone by one unit |
| $Q_c$ | Delivery target at community $c \in C$ |
| $L_{max}$ | Maximum delivery capacity of a delivery drone |
| $\bar{x}_{ij}$ | =1, if a truck traverses arc $(i,j) \in E^o$ in the final solution of the first stage; 0 otherwise |
| $\mathrm{T}_{sc}^d$ | Time to reach community $c \in C$ by a delivery drone launched from satellite $s \in S$ |
| $e_{sc}$ | =1, if community $c \in C$ is provided telecommunication coverage by a hotspot drone launched from satellite $s \in S$ |
| $\bar{\Delta}_s^h$ | The longest route completion time by a hotspot drone launched from satellite $s \in S$ in the first stage |

$d_{sk}^{c}$ Delivered quantity at community $c \in C$ by a delivery drone following route $k \in V_s'$ from satellite $s \in S$

*Newly defined decision variables*

$A_c$ Missed delivery amount at community $c \in C$

$G_{sk}$ Carried load amount exceeding the drone load capacity by delivery drone $k \in D_s$ launched from satellite $s \in S$

We keep the variables $z_l^{sk}, \Delta_s^d, T_s, T_{sc}^d;\ \forall s \in S, l \in V_s', k \in D_s, c \in C$ defined in the first stage. We keep the constraints (5), (7), (9), (12), (13) from the first stage model with the modification below:

fix $x_{ij} \coloneqq \bar{x}_{ij};\ \forall (i,j) \in E^o$ and $\Delta_s^h := \bar{\Delta}_s^h;\ \forall s \in S$

The restricted main problem for the second stage (D-RMP-2) can now be presented in (23)—(28).

$$Minimize \sum_{s\in S}\sum_{c\in C} F_c^T \mathrm{T}_{sc}^d + \sum_{c\in C} F_c^R J_c + \sum_{c\in C} F_c^D A_c + \sum_{s\in S}\sum_{k\in D_s} F^E G_{sk} \quad (23)$$

Subject to: **Duals**

$$(5), (7), (9), (12), (13) \quad (24)$$

$$\sum_{s\in S}\sum_{k\in D_s}\sum_{l\in V_s'} z_l^{sk}\, d_{sk}^c + A_c \geq Q_c \qquad \forall c \in C \qquad \zeta_c' \quad (25)$$

$$\sum_{k\in D_s}\sum_{l\in V_s'} z_l^{sk}\, b_c^l - e_{sc} \leq 0 \qquad \forall s \in S, c \in C, \qquad \pi_{sc}' \quad (26)$$

$$\sum_{l\in V_s'}\sum_{c\in C} z_l^{sk} d_{sk}^c \leq L_{max} + G_{sk} \qquad \forall s \in S, k \in D_s \qquad \mu_{sk}' \quad (27)$$

$$z_l^{sk} \in \{0,1\};\ \Delta_s^d, T_s, \mathrm{T}_{sc}^d, G_{sk}, A_c \geq 0$$
$$\forall s \in S, l \in V_s', k \in D_s, c \in C \quad (28)$$

The second stage objective function (23) minimizes: (i) the cost associated with the time to reach communities by delivery drones, (ii) the penalty for failure to reach communities by delivery drones, (iii) the cost of unfulfillment of demand at communities, and (iv) the cost of violating delivery drone load carrying capacities. Constraints (25) ensure that the demand at community $c \in C$ must be satisfied by the delivery drone or variable $A_c$ becomes equal to the value of the missed delivery amount. According to constraints (26), community $c \in C$ can be served by a delivery drone launched from satellite $s \in S$ only if community $c$ is provided telecommunication coverage by a hotspot drone launched from $s$. Constraints (27) ensures that the total delivery amount by drone $k \in D_s$; $\forall s \in S$ at the communities along its route is either bounded by the drone load capacity, or variable $G_{sk}$ takes the value equal to the excess amount. Lastly, the variable restrictions are given in (28). Like in the first stage, the generation of new delivery routes and determination of route parameters for delivery drones are relegated to the pricing problem.

### 3.5. Pricing subproblems for generating delivery drone routes-2nd stage (D-PSP-DD-2)

We introduce two new sets of decision variables for the pricing problems.

*Newly defined decision variables:*

$q_{ij}^{sk}$ Quantity transported through arc $(i,j) \in E_s^d$ by a delivery drone launched from satellite $s \in S$ by delivery drone $k \in D_s$

$p_j^{sk}$ Amount delivered at community $j \in C$ by a delivery drone launched from satellite $s \in S$ by delivery drone $k \in D_s$

The pricing subproblems for generating route for each delivery drone $k \in D_s$ from satellite $s \in S$, are given in (29)—(34).

$$\begin{aligned} Minimize \quad (RC')_{sk}^d &= -\theta_s' - \phi_{sk}' - \sum_{j \in C} \sum_{(i,j) \in E_s^d} v_{ij}^{sk} \left(\pi_{sj}' + \rho_{sk}' \tau_{ij}^d\right) \\ &\quad - \sum_{j \in C_s} \alpha_{sj}' \left( \mathrm{t}_j^{sk} + M \sum_{(i,j) \in E_s^d} v_{ij}^{sk} \right) - \sum_{j \in C_s} (\zeta_j' + \mu_{sk}') p_j^{sk} \end{aligned} \tag{29}$$

Subject to:

$$(17)-(22) \tag{30}$$

$$q_{ij}^{sk} \leq v_{ij}^{sk} L_{max} \qquad \forall\, (i,j) \in E_s^d, \tag{31}$$

$$\sum_{(i,j)\in E_s^d} q_{ij}^{sk} - \sum_{(j,k)\in E_s^d} q_{jk}^{sk} = p_j^{sk} \qquad \forall\, j \in C_s, \tag{32}$$

$$\sum_{j\in C} p_j^{sk} \leq L_{max} \tag{33}$$

$$p_j^{sk}, q_{ij}^{sk} \geq 0 \qquad \forall j \in C_s, (i,j) \in E_s^d \tag{34}$$

Objective function (29) minimizes the reduce cost $(RC')_{sk}^{d}$ of a delivery route for drone $k \in D_s$ launched from satellite $s \in S$, where $\zeta_c', \pi_{sc}', \mu_{sk}'$ are the dual variables associated with constraints (25)—(27), and $\theta_s', \phi_{sk}', \rho_{sk}', \alpha_{sc}'$ are the dual variables associated with constraints (7), (9), (12), (13) as modified for D-RMP-2. If the reduced cost is found to be negative, then the route is added to the set $V_s'$. Constraints (31)—(33) ensure that load carrying capacity of delivery drones are not violated and restrict the delivery amount at each community. In the next subsection, we present our solution algorithm for this two-stage decision framework.

### *3.6. Algorithm for solving the deterministic two-echelon vehicle routing model*

**Step 1.** We initialize an iteration counter, generate initial routes for the hotspot and delivery drones and add them to $U_s; \forall s \in S$ and $V_s; \forall s \in S$, respectively. These initial sets can be empty or can contain routes with only one hotspot node and one community node, respectively.

**Step 2.** We relax the variable integrality restrictions and solve the D-RMP-1 with the current restricted route sets $U_s, V_s; \forall s \in S$. We obtain dual solution sets $\gamma_s$, $\theta_s$, $\varphi_{sk}$, $\phi_{sl}$, $\zeta_c$, $\sigma_{sk}$, $\rho_{sl}$, $\alpha_{sc}$, $\pi_{sc}; \forall s \in S, k \in H_s, l \in D_s, c \in C$ associated with constraints (6)—(14).

**Step 3a.** We solve D-PSP-HD-1 (see Appendix) for each satellite and each hotspot drone ($\forall s \in S, k \in H_s$). Following each pricing problem solution, we update the hotspot node set $B_s \coloneqq B_s -$

$\left\{j \in B_s : \sum_{(i,j)\in E_s^h} u_{ij}^{s(k-1)} = 1\right\}; \forall k \in H_s : k = 2., \dots, m^h$. For each satellite, we continue solving D-PSP-DD-1 until $B_s = \emptyset$, or if number of routes generated equals $m^h$. If $\exists s \in S, k \in H_s : (RC)_{sk}^h \geq 0$, we store the corresponding hotspot drone route to $U_s; \forall s \in S$.

**Step 3b.** We solve D-PSP-DD-1 for each satellite and each delivery drone ($\forall s \in S, k \in D_s$). Following each pricing problem solution, we update the community node set $C_s \coloneqq C_s - \left\{j \in C_s : \sum_{(i,j)\in E_s^d} v_{ij}^{s(k-1)} = 1\right\}; \forall k \in D_s : k = 2., \dots, m^d$. For each satellite, we continue solving D-PSP-DD-2 until $C_s = \emptyset$, or if number of routes generated equals $m^d$. If $\exists s \in S, k \in D_s : (RC)_{sk}^d \geq 0$, we store the corresponding delivery drone route to $V_s; \forall s \in S$.

**Step 4.** If $(RC)_{sk}^h \geq 0;\ \forall s \in S, k \in H_s$ in *Step 3a,* and $(RC)_{sk}^d \geq 0;\ \forall s \in S, k \in D_s$ in *Step 3b*, go to *Step 5.* Otherwise, if $\exists s \in S, k \in H_s : (RC)_{sk}^h < 0$ in *Step 3a*, or $\exists s \in S, k \in D_s : (RC)_{sk}^d < 0$ in *Step 3b*, we go back to *Step 2.*

**Step 5.** We solve the D-RMP-1 with the current restricted sets $U_s, V_s; \forall s \in S$ using the relax-and-fix heuristic. At first, D-RMP-1 is solved with relaxing the integrality restriction on variables for delivery drone routes ($\boldsymbol{z}$) and the optimal solution for truck routes ($\boldsymbol{x}$) and hotspot drone routes ($\boldsymbol{y}$) are obtained. Next, the variables $\boldsymbol{x}$ and $\boldsymbol{y}$ are fixed at their optimal values and the D-RMP-1 is solved again with integrality restriction imposed on $\boldsymbol{z}$. We store the solutions from D-RMP-1 for the second stage as follows: $\bar{x}_{ij} \coloneqq x_{ij}^*;\ \forall (i,j) \in E^o;\ V_s' \coloneqq V_s^*, \forall s \in S;\ \bar{\Delta}_s^h \coloneqq \Delta^{*h}_s, \forall s \in S;\ e_{sc} = \min\left(1, \sum_{l\in U_s} \sum_{k\in H_s} \sum_{j\in B_s} a_j^l\, y_l^{sk} g_{jc}\right), \forall s \in S, c \in C,\ d_{sk}^c := b_c^k Q_c, \forall s \in S, k \in V_s', c \in C$ and then go to *Step 6.*

**Step 6.** We relax the variable integrality restrictions and solve the D-RMP-2 with the current restricted route sets $V_s'; \forall s \in S$. We obtain dual solutions $\zeta_c', \pi_{sc}', \mu_{sk}'$ associated with constraints (25)—(27), and also $\theta_s'$, $\phi_{sk}'$, $\rho_{sk}', \alpha_{sc}'$ that are associated with constraints (7), (9), (12), (13) modified for D-RMP-2.

**Step 7.** We solve D-PSP-DD-2 for each satellite and each delivery drone ($\forall s \in S, k \in D_s$). Following each pricing problem solution, we update the community node set the same way as we did in *Step 3b*. For each satellite, we continue solving D-PSP-DD-2 until $C_s = \emptyset$, or if number of routes generated equals $m^d$. If $\exists s \in S, k \in D_s : (RC)_{sk}^h \geq 0$, we store the corresponding delivery drone route to $V_s'$, $\forall s \in S$.

**Step 8.** If $(RC')_{sk}^d \geq 0, \forall s \in S, k \in D_s$ in *Step 7*, go to *Step 9*. Otherwise, if $\exists s \in S, k \in D_s$ : $(RC')_{sk}^d < 0$, we go back to *Step 6*.

**Step 9.** Imposing variable integrality restrictions, we solve the D-RMP-2 with the current restricted drone route set $V_s', \forall s \in S$. We stop the algorithm when the optimal delivery drone route set is found by D-RMP-2.

## 4. Numerical experiments and discussion

In this section, we present a case study for post-disaster emergency aid deliveries in the island of Puerto Rico using our framework. We visually present the selection of satellite, hotspot, and demand locations in Figure 2. Satellite stations, corresponding to the intersections of major highways, are shown in red, hotspot locations, corresponding to pre-disaster telecommunication towers, are in green, and delivery locations, corresponding to zip codes, are in blue.

Two of the most important metrics for post-disaster aid delivery operations are the **average proportion of unfulfilled demand** and the **average delay in reaching communities**. In the remainder of this section, we investigate the effect of various parameters on these two metrics using a case study simulating a post-disaster emergency aid delivery problem in Puerto Rico.

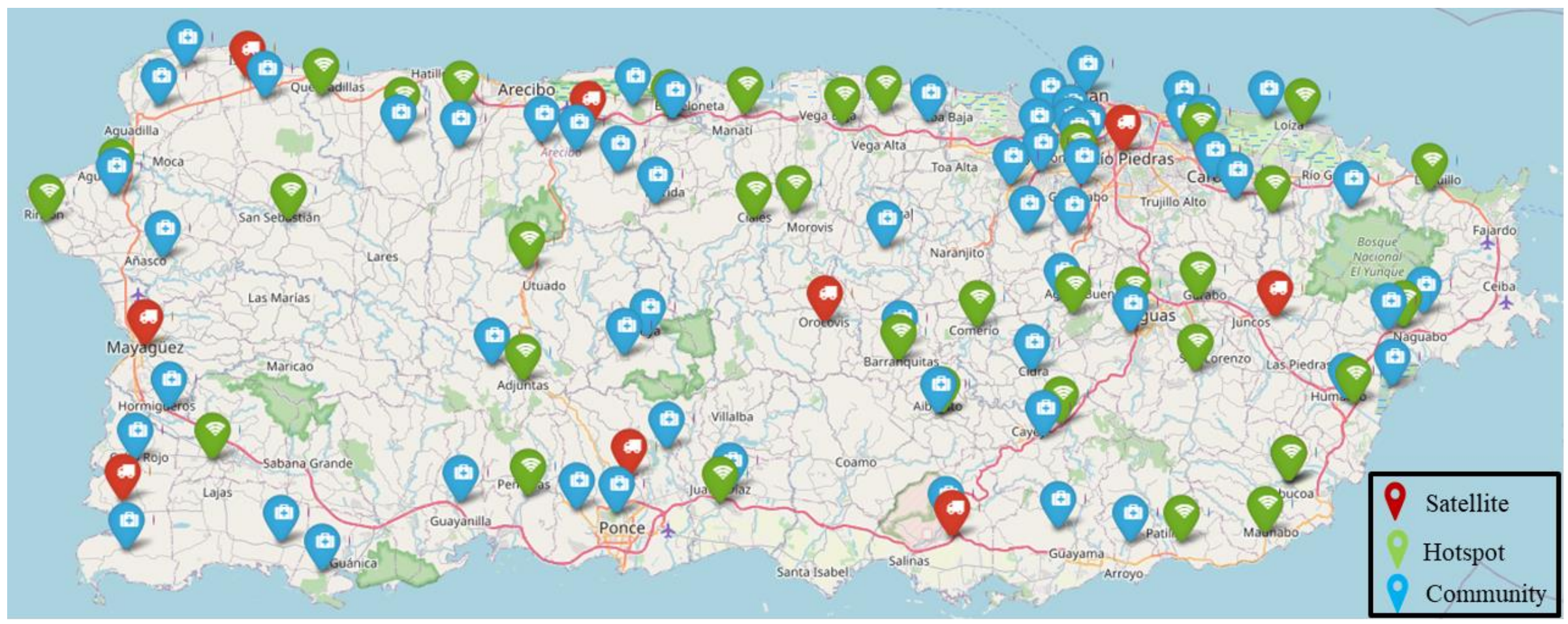

Figure 2: Locations of satellites, hotspots, and communities.

Some of the parameters are to be selected by policy makers and stakeholders; others are dependent on the nature and severity of the disaster and the geography of the impacted area. Our numerical experiments and result analyses can provide humanitarian logistics service providers with crucial insights for the necessary trade-offs between unfulfilled demand and delay in sending emergency aid.

### *4.1. Computational experiments*

We create datasets of four different sizes, with 60, 80, 100, and 120 community locations, 10 satellite locations, and 40 hotspot locations. For each problem instance, we randomly select the locations of communities and generate demand quantities in the range of $[2,15]$. We set the number of trucks at 4 (i.e., $m^t = 4$) and consider 1 hotspot drone in each truck ($m^h = 1$). For choosing the number of delivery drones per truck ($m^d$), first we fix the delivery drone load capacity, $L_{max}$, at 25 units, which exceeds the maximum possible demand quantity at any community. Next, we identify the minimum number of drones that are necessary to exceed the sum of demand quantities at the communities that can be served from each satellite, and we select the minimum of these numbers to serve as the number of delivery drones carried by a truck. We set the drone flying ranges $W^h$ and $W^d$ at 35 miles. To calculate the distances between nodes, we use the geodesic distance measure and consider 1 mile/minute for all ground and aerial vehicles. For binary parameters $g_{jc}, \forall j \in B, c \in C$, i.e., we consider three different hotspot coverage radii (signal strengths) for hotspot drones equal to 12, 16, and 20 miles. We set delay cost parameters $F_c^T$ at \$1/minute, $F_c^R$ at \$10000 for failure to reach a community, and $F^E$ at \$10000 for per unit violation of delivery drone load capacity. For the per unit demand unfulfillment cost $F_c^D$ parameter, we randomly generate it in the range of $[10, 1000]$.

We implement the model and solution algorithm in AMPL and then use the CPLEX 20.1.0 solver to obtain optimal solutions to all generated datasets. We use a computer with Dual Intel Xeon Processor (12 Core, 2.3GHz Turbo) with 64 GB of RAM, and Windows 64-bit operating system. A summary of our results is presented in Table 1. To run these experiments, we calculate $m^d$ using the approach mentioned above and set its value at 4, 6, 9, and 10, for 60, 80, 100, and 120 size instances respectively.

From Table 1, we see that the total cost decreases with the increase in hotspot coverage radius for each instance size, as expected. We further notice that when the hotspot coverage radius is small and fewer delivery drones are carried in each truck, a higher proportion of the total requested demand is unfulfilled. For example, when the coverage radius is 12 miles, the 60-community case with only 4 drones in each truck results in a higher proportion of unfulfilled demand compared to larger size instances, when more drones are carried in a truck. An interesting outcome from our experimentation is that as the number of drones per truck increases, the impact of the hotspot coverage radius diminishes. We also notice the effect of hotspot coverage radius on the average delay to reach communities. A larger coverage radius results in a more relaxed version of the problem and usually leads to decreased average delay times. We consider delay time as the length of time span between the first truck staring its route (the beginning of the planning horizon) and a delivery drone reaching a community.

Table 1: Computational results summary for the 2EVRP-HD-UAV framework (5 case average).

| Instance size | Hotspot coverage radius (miles) | CPU time (second) | | Total cost of delay and unfulfillment ($) | | Proportion of unfulfilled demand (%) | | Average delay to reach communities (minute) | |
|---|---|---|---|---|---|---|---|---|---|
| | | Mean | Standard deviation | Mean | Standard deviation | Mean | Standard deviation | Mean | Standard deviation |
| 60 | 12 | 399.40 | 32.02 | 13804.60 | 5188.57 | 11.05 | 4.33 | 89.85 | 2.55 |
| | 16 | 241.20 | 30.56 | 5765.80 | 324.61 | 0.00 | 0.00 | 96.10 | 5.41 |
| | 20 | 255.00 | 18.03 | 5360.20 | 336.22 | 0.00 | 0.00 | 89.34 | 5.60 |
| 80 | 12 | 1428.60 | 393.38 | 14746.60 | 2050.49 | 7.66 | 2.00 | 94.13 | 0.92 |
| | 16 | 369.60 | 32.03 | 7684.00 | 178.77 | 0.00 | 0.00 | 96.05 | 2.23 |
| | 20 | 538.25 | 38.40 | 7502.40 | 141.51 | 0.00 | 0.00 | 93.78 | 1.77 |
| 100 | 12 | 2255.00 | 453.32 | 16644.60 | 6106.53 | 6.99 | 4.14 | 91.53 | 0.47 |
| | 16 | 1072.40 | 80.00 | 8786.60 | 213.64 | 0.00 | 0.00 | 87.87 | 2.14 |
| | 20 | 1200.20 | 87.34 | 8429.00 | 165.20 | 0.00 | 0.00 | 84.29 | 1.65 |
| 120 | 12 | 3183.00 | 123.37 | 12622.60 | 289.01 | 0.43 | 0.10 | 104.23 | 2.62 |
| | 16 | 2359.00 | 264.59 | 10518.00 | 117.35 | 0.00 | 0.00 | 87.65 | 1.65 |
| | 20 | 2939.60 | 170.77 | 10100.40 | 162.34 | 0.00 | 0.00 | 84.17 | 1.35 |

### *4.2 Sensitivity analysis*

We run experiments to study the effect of various parameters on the routing decisions and the resultant cost, the proportion of unfulfilled demand, and the average delay in reaching communities. We focus on the 60-community size instances for this analysis and present five datasets. To generate the first four datasets, we randomly select communities across Puerto Rico with higher probability to different geographical areas: (i) northwestern part of the island, (ii) northeastern part of the island, (iii) southwestern part of the island, and (iv) southeastern part of the island, respectively. We present these datasets in the maps of Figure 3. For the fifth dataset, we randomly select locations without any geographical preference (i.e., all locations are equally probably selected).

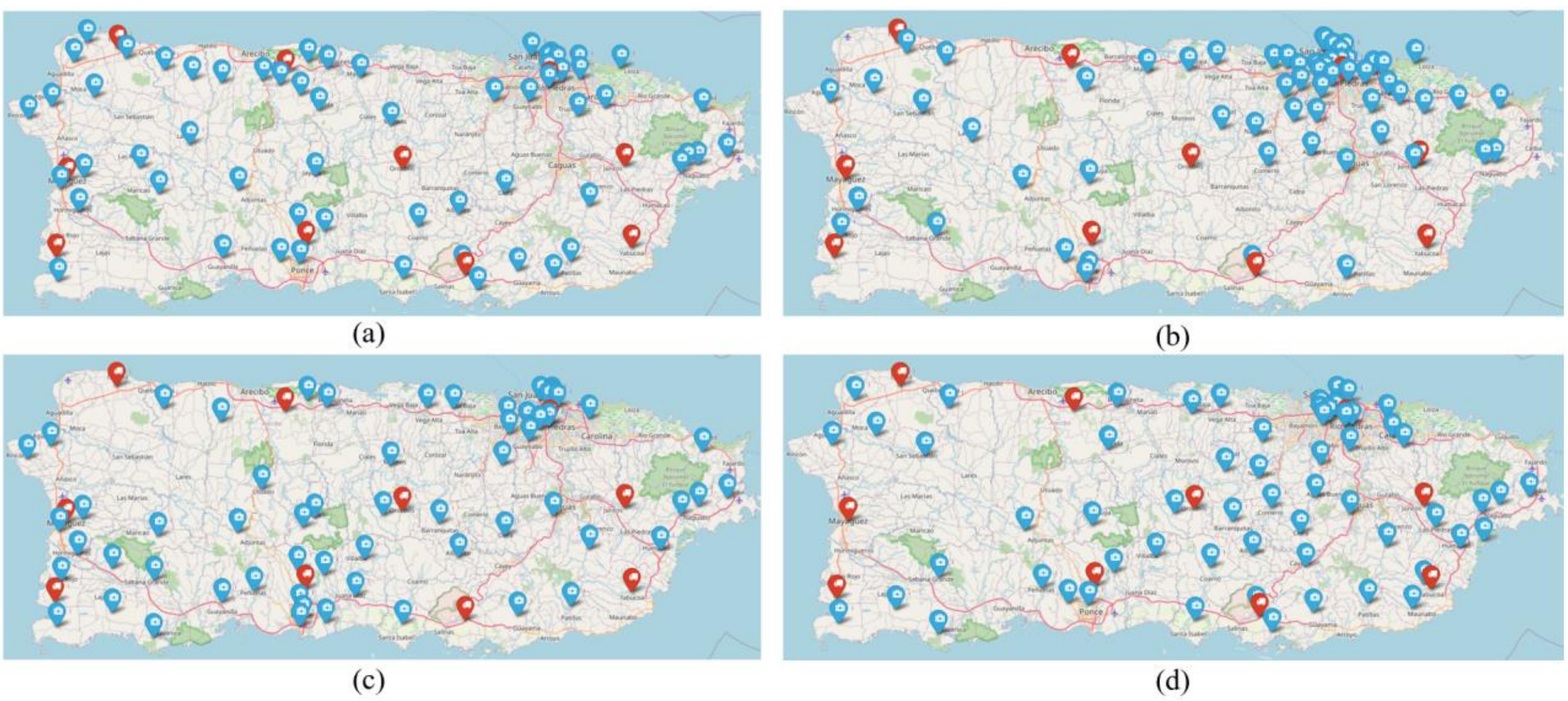

Figure 3: Locations of satellites (red pins), and communities (blue pins) predominantly in (a) northwestern, (b) northeastern, (c) southwestern, and (d) southeastern regions of Puerto Rico.

We run numerical experiments using these datasets. For the 2EVRP-HD-UAV framework, we introduce several variations of the number of delivery drones carried in each truck and flying ranges of delivery drones. We find that the unfulfilled proportion of demand decreases with the increase in the delivery drone flying range and the number of delivery drones carried in a truck. The rate of decrease varies with the hotspot coverage radius, and the unfulfilled proportion of demand decreases at a highest rate when the hotspot coverage radius is set at 12 miles. For higher values of coverage radius, e.g., for 16 and 20 miles, the rate of decrease is very similar. We present these insights in Figure 4.

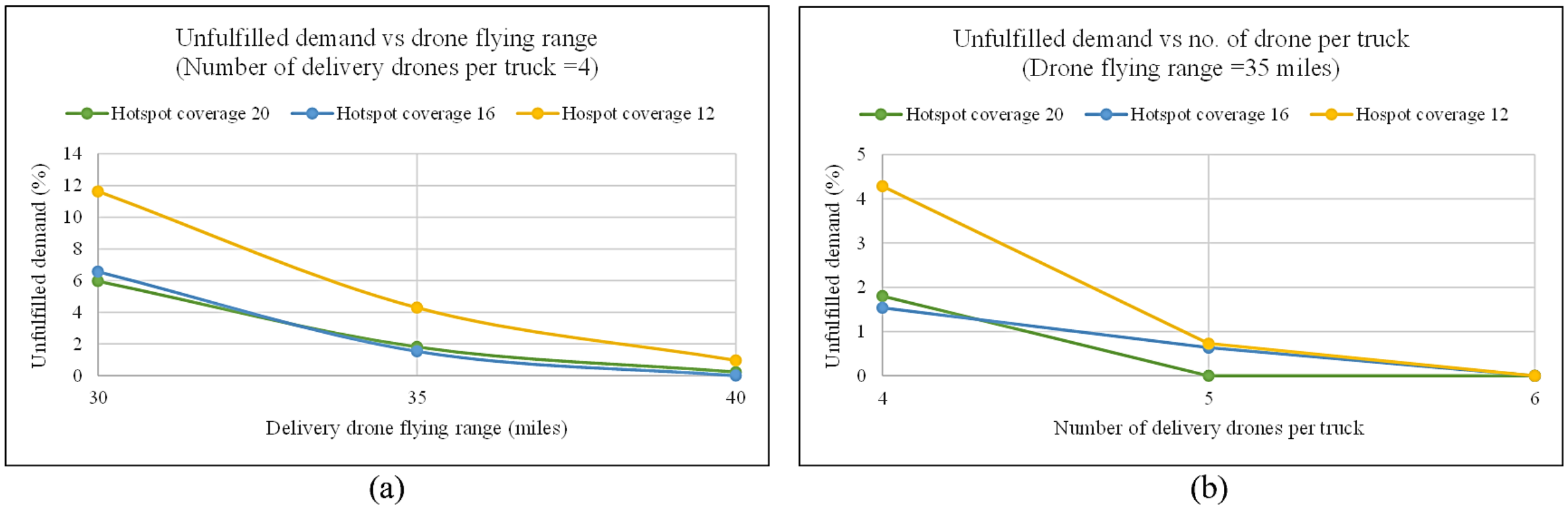

Figure 4: Variation in unfulfilled demand with change in (a) delivery drone flying range and (b) number of delivery drones per truck.

We also study the variation in the average delay in reaching communities as the delivery drone flying range and the number of delivery drones carried in a truck change. For a fixed number of delivery drones carried by a truck, the average delay time increases with the increase in drone flying range, whereas the proportion of unfulfilled demand decreases as evident from Figure 5(a). The increase in the delay time halts if the unfulfilled proportion of demand reaches zero. On the other hand, both the average delay and the unfulfilled proportion of demand decrease with an increase in the number of delivery drones, which can be seen in Figure 5(b).

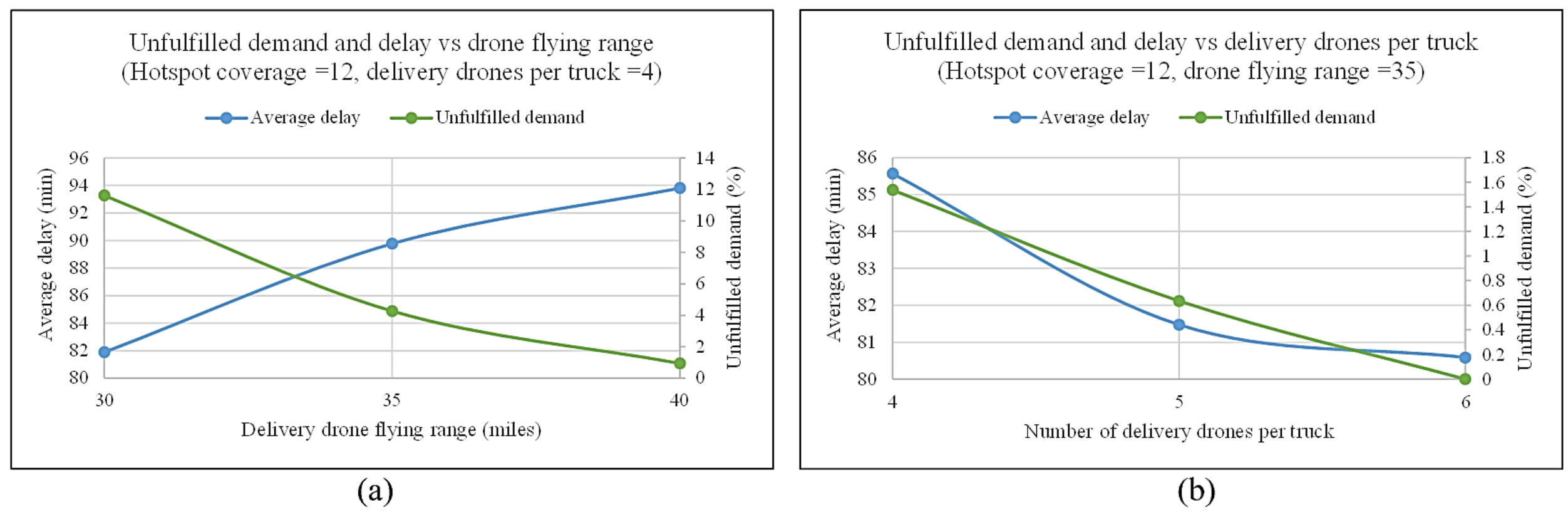

Figure 5: Variation in unfulfilled demand and average delay in reaching communities with change in (a) delivery drone flying range and (b) number of delivery drones per truck.

## 5. Concluding remarks

In this paper, we study two-echelon vehicle routing problems (2E-VRP) under stochastic conditions in the context of humanitarian logistics. Following a disaster, the demand for emergency supplies (medical kits, dry food, water) in the affected communities is uncertain due to absence of communication infrastructure (to communicate the demands); furthermore, the delivery of the necessary supplies through traditional ground transportation is not always possible due to infrastructure damage. To handle this humanitarian aid delivery problem, we propose a unique 2E-VRP framework, where trucks carrying drones travel from satellite location to satellite location. When a truck stops at one of those designated locations, drones are dispatched for second echelon ("last mile") deliveries. We develop a computational framework to handle the demand uncertainty, where a special purpose drones, i.e., hotspot drones are used for handling demand uncertainty, by providing communications capabilities to allow people to communicate their demands.

In the presented framework, we propose a solution approach in two stages. In the first stage, we identify optimal routes for reaching every demand location such that total time to reach a community is minimized, while we also cover each community with the use of a hotspot drone (for providing a communication signal). In the second stage, the demand information is revealed, and we choose among the routes obtained in the first stage or identify better routes to minimize the cost of reaching the demand locations as well as the cost of missing demands. We design a decomposition scheme and an efficient column generation-based heuristic solution approach that results in small and manageable pricing subproblems for hotspot and delivery drone route generation.

In our current 2E-VRP framework, we consider only one source of uncertainty, i.e., demand for emergency aid at the affected communities. During real-world humanitarian logistics operations, several other parameters are of a stochastic nature. For example, we consider the flying time range and load carrying capacities of drones as deterministic. These two properties of drones are significantly impacted by weather conditions during delivery operations. We plan to extend our study by incorporating uncertain weather conditions in terms of stochastic flying range and load carrying capacity of drones in the future.

## Conflict of interest

We have no conflict of interest to declare.

## References

[1] G. Kovács and K. M. Spens, "Humanitarian logistics in disaster relief operations, "*International Journal of Physical Distribution & Logistics Management*, vol. 37, no. 2, pp. 99-114, 2007.

[2] A. Leiras, I. de Brito Jr, E. Queiroz Peres, T. Rejane Bertazzo, and H. Tsugunobu Yoshida Yoshizaki, "Literature review of humanitarian logistics research: trends and challenges," *Journal of Humanitarian Logistics and Supply Chain Management*, vol. 4, no. 1, pp. 95-130, 2014.

[3] R. Banomyong, P. Varadejsatitwong, and R. Oloruntoba, "A systematic review of humanitarian operations, humanitarian logistics and humanitarian supply chain performance literature 2005 to 2016," *Annals of Operations Research*, vol. 283, no. 1-2, pp. 71-86, 2019.

[4] L. Özdamar, E. Ekinci, and B. Küçükyazici, "Emergency logistics planning in natural disasters," *Annals of operations research,* vol. 129, no. 1-4, pp. 217-245, 2004.

[5] L. Özdamar and O. Demir, "A hierarchical clustering and routing procedure for large scale disaster relief logistics planning," *Transportation Research Part E: Logistics and Transportation Review,* vol. 48, no. 3, pp. 591-602, 2012.

[6] H. Wang, L. Du, and S. Ma, "Multi-objective open location-routing model with split delivery for optimized relief distribution in post-earthquake," *Transportation Research Part E: Logistics and Transportation Review,* vol. 69, pp. 160-179, 2014.

[7] C. C. Lu, K. C. Ying, and H. J. Chen, "Real-time relief distribution in the aftermath of disasters–A rolling horizon approach," *Transportation research part E: logistics and transportation review,* vol. 93, pp. 1-20, 2016.

[8] T. I. Faiz, C. Vogiatzis, and M. Noor-E-Alam, "A column generation algorithm for vehicle scheduling and routing problems," *Computers & Industrial Engineering,* vol. 130, pp. 222-236, 2019.

[9] R. Cuda, G. Guastaroba, and M. G. Speranza, "A survey on two-echelon routing problems," *Computers & Operations Research,* vol. 55, pp. 185-199, 2015/03/01/ 2015, doi: https://doi.org/10.1016/j.cor.2014.06.008.

[10] J. Gonzalez-Feliu, G. Perboli, R. Tadei, and D. Vigo, "The two-echelon capacitated vehicle routing problem," 2008. halshs-00879447.

[11] G. Perboli and R. Tadei, "New families of valid inequalities for the two-echelon vehicle routing problem," *Electronic notes in discrete mathematics,* vol. 36, pp. 639-646, 2010.

[12] G. Perboli, R. Tadei, and D. Vigo, "The two-echelon capacitated vehicle routing problem: models and math-based heuristics," *Transportation Science,* vol. 45, no. 3, pp. 364-380, 2011.

[13] M. Jepsen, S. Spoorendonk, and S. Ropke, "A branch-and-cut algorithm for the symmetric two-echelon capacitated vehicle routing problem," *Transportation Science,* vol. 47, no. 1, pp. 23-37, 2013.

[14] R. Baldacci, A. Mingozzi, R. Roberti, and R. W. Calvo, "An exact algorithm for the two-echelon capacitated vehicle routing problem," *Operations Research,* vol. 61, no. 2, pp. 298-314, 2013.

[15] F. A. Santos, G. R. Mateus, and A. S. da Cunha, "A branch-and-cut-and-price algorithm for the two-echelon capacitated vehicle routing problem," *Transportation Science,* vol. 49, no. 2, pp. 355-368, 2014.

[16] F. A. Santos, A. S. da Cunha, and G. R. Mateus, "Branch-and-price algorithms for the Two-Echelon Capacitated Vehicle Routing Problem," *Optimization Letters,* journal article vol. 7, no. 7, pp. 1537-1547, October 01 2013, doi: 10.1007/s11590-012-0568-3.

[17] T. Liu, Z. Luo, H. Qin, and A. Lim, "A branch-and-cut algorithm for the two-echelon capacitated vehicle routing problem with grouping constraints," *European Journal of Operational Research,* vol. 266, no. 2, pp. 487-497, 2018/04/16/ 2018, doi: https://doi.org/10.1016/j.ejor.2017.10.017.

[18] G. Perboli, R. Tadei, and E. Fadda, "New Valid Inequalities for the Two-Echelon Capacitated Vehicle Routing Problem," *Electronic Notes in Discrete Mathematics,* vol. 64, pp. 75-84, 2018/02/01/ 2018, doi: https://doi.org/10.1016/j.endm.2018.01.009.

[19] T. G. Crainic, G. Perboli, S. Mancini, and R. Tadei, "Two-echelon vehicle routing problem: a satellite location analysis," *Procedia-Social and Behavioral Sciences,* vol. 2, no. 3, pp. 5944-5955, 2010.

[20] T. G. Crainic, S. Mancini, G. Perboli, and R. Tadei, "Multi-start heuristics for the two-echelon vehicle routing problem," in *European Conference on Evolutionary Computation in Combinatorial Optimization*, 2011: Springer, pp. 179-190.

[21] V. C. Hemmelmayr, J.-F. Cordeau, and T. G. Crainic, "An adaptive large neighborhood search heuristic for two-echelon vehicle routing problems arising in city logistics," *Computers & operations research,* vol. 39, no. 12, pp. 3215-3228, 2012.

[22] U. Breunig, V. Schmid, R. F. Hartl, and T. Vidal, "A large neighbourhood based heuristic for two-echelon routing problems," *Computers & Operations Research,* vol. 76, pp. 208-225, 2016/12/01/ 2016, doi: https://doi.org/10.1016/j.cor.2016.06.014.

[23] L. Zhou, R. Baldacci, D. Vigo, and X. Wang, "A multi-depot two-echelon vehicle routing problem with delivery options arising in the last mile distribution," *European Journal of Operational Research,* vol. 265, no. 2, pp. 765-778, 2018.

[24] Z.-y. Zeng, W.-s. Xu, Z.-y. Xu, and W.-h. Shao, "A Hybrid GRASP+VND Heuristic for the Two-Echelon Vehicle Routing Problem Arising in City Logistics," *Mathematical Problems in Engineering,* vol. 2014, p. 11, 2014, Art no. 517467, doi: 10.1155/2014/517467.

[25] M. Drexl, "Synchronization in Vehicle Routing—A Survey of VRPs with Multiple Synchronization Constraints," *Transportation Science,* vol. 46, no. 3, pp. 297-316, 2012, doi: 10.1287/trsc.1110.0400.

[26] S. Martinez, "UAV Cooperative Decision and Control: Challenges and Practical Approaches (Shima, T. and Rasmussen, S.; 2008) [Bookshelf]," *IEEE Control Systems Magazine,* vol. 30, no. 2, pp. 104-107, 2010, doi: 10.1109/MCS.2010.935899.

[27] S. G. Manyam, D. W. Casbeer, and K. Sundar, "Path planning for cooperative routing of air-ground vehicles," in *2016 American Control Conference (ACC)*, 6-8 July 2016 2016, pp. 4630-4635, doi: 10.1109/ACC.2016.7526082.

[28] Z. Luo, Z. Liu, and J. Shi, "A Two-Echelon Cooperated Routing Problem for a Ground Vehicle and Its Carried Unmanned Aerial Vehicle," *Sensors (Basel, Switzerland),* vol. 17, no. 5, p. 1144, 2017, doi: 10.3390/s17051144.

[29] X. Wang, S. Poikonen, and B. Golden, "The vehicle routing problem with drones: several worst-case results," *Optimization Letters,* vol. 11, no. 4, pp. 679-697, 2017/04/01 2017, doi: 10.1007/s11590-016-1035-3.

[30] S. Poikonen, X. Wang, and B. Golden, "The vehicle routing problem with drones: Extended models and connections," *Networks,* vol. 70, no. 1, pp. 34-43, 2017, doi: doi:10.1002/net.21746.

[31] B. Behdani and J. C. Smith, "An Integer-Programming-Based Approach to the Close-Enough Traveling Salesman Problem," *INFORMS Journal on Computing,* vol. 26, no. 3, pp. 415-432, 2014, doi: 10.1287/ijoc.2013.0574.

[32] J. G. Carlsson and S. Song, "Coordinated Logistics with a Truck and a Drone," *Management Science,* vol. 64, no. 9, pp. 4052-4069, 2018, doi: 10.1287/mnsc.2017.2824.

[33] C. C. Murray and A. G. Chu, "The flying sidekick traveling salesman problem: Optimization of drone-assisted parcel delivery," *Transportation Research Part C: Emerging Technologies,* vol. 54, pp. 86-109, 2015/05/01/ 2015, doi: https://doi.org/10.1016/j.trc.2015.03.005.

[34] N. Agatz, P. Bouman, and M. Schmidt, "Optimization Approaches for the Traveling Salesman Problem with Drone," *Transportation Science,* vol. 52, no. 4, pp. 965-981, 2018, doi: 10.1287/trsc.2017.0791.

[35] E. Es Yurek and H. C. Ozmutlu, "A decomposition-based iterative optimization algorithm for traveling salesman problem with drone," *Transportation Research Part C: Emerging Technologies,* vol. 91, pp. 249-262, 2018/06/01/ 2018, doi: https://doi.org/10.1016/j.trc.2018.04.009.

[36] Q. M. Ha, Y. Deville, Q. D. Pham, and M. H. Ha, "On the min-cost traveling salesman problem with drone," *Transportation Research Part C: Emerging Technologies,* vol. 86, pp. 597-621, 2018.

[37] Y. S. Chang and H. J. Lee, "Optimal delivery routing with wider drone-delivery areas along a shorter truck-route," *Expert Systems with Applications,* vol. 104, pp. 307-317, 2018/08/15/ 2018, doi: https://doi.org/10.1016/j.eswa.2018.03.032.

[38] P. Kitjacharoenchai, M. Ventresca, M. Moshref-Javadi, S. Lee, J. M. A. Tanchoco, and P. A. Brunese, "Multiple traveling salesman problem with drones: Mathematical model and heuristic approach," *Computers & Industrial Engineering,* vol. 129, pp. 14-30, 2019/03/01/ 2019, doi: https://doi.org/10.1016/j.cie.2019.01.020.

[39] M. Hu *et al.*, "Joint Routing and Scheduling for Vehicle-Assisted Multi-Drone Surveillance," *IEEE Internet of Things Journal,* pp. 1-1, 2019, doi: 10.1109/JIOT.2018.2878602.

[40] H. Savuran and M. Karakaya, *Route Optimization Method for Unmanned Air Vehicle Launched from a Carrier*. 2015, pp. 279-284.

[41] H. Savuran and M. Karakaya, "Efficient route planning for an unmanned air vehicle deployed on a moving carrier," *Soft Computing,* journal article vol. 20, no. 7, pp. 2905-2920, July 01 2016, doi: 10.1007/s00500-015-1970-4.

[42] Z. Ghelichi, M. Gentili, and P. B. Mirchandani, "Logistics for a fleet of drones for medical item delivery: A case study for Louisville, KY," *Computers & Operations Research,* p. 105443, 2021.

[43] L. Evers, T. Dollevoet, A. I. Barros, and H. Monsuur, "Robust UAV mission planning," *Annals of Operations Research,* vol. 222, no. 1, pp. 293-315, 2014/11/01 2014, doi: 10.1007/s10479-012-1261-8.

[44] S. J. Kim, G. J. Lim, and J. Cho, "Drone flight scheduling under uncertainty on battery duration and air temperature," *Computers & Industrial Engineering,* vol. 117, pp. 291-302, 2018/03/01/ 2018, doi: https://doi.org/10.1016/j.cie.2018.02.005.

[45] K. Wang, S. Lan, and Y. Zhao, "A genetic-algorithm-based approach to the two-echelon capacitated vehicle routing problem with stochastic demands in logistics service," *Journal of the Operational Research Society,* vol. 68, no. 11, pp. 1409-1421, 2017.

[46] R. Liu, Y. Tao, Q. Hu, and X. Xie, "Simulation-based optimisation approach for the stochastic two-echelon logistics problem," *International Journal of Production Research,* vol. 55, no. 1, pp. 187-201, 2017.

## Appendix

### *A1. Mixed integer linear programming model*

The full MILP model is presented in (A1)—(A24).

$$Minimize \sum_{s\in S}\sum_{c\in C} F_c^T \mathrm{T}_{sc}^d + \sum_{c\in C} F_c^R J_c + \sum_{c\in C} F_c^D A_c \qquad \text{(A1)}$$

Subject to:

$$\sum_{(i,j)\in E^o: i\notin S} x_{ij} \leq m^t \qquad \text{(A2)}$$

$$\sum_{(i,j)\in E^o} x_{ij} \leq 1 \qquad \forall j \in S, \qquad \text{(A3)}$$

$$\sum_{(i,j)\in E^o} x_{ij} - \sum_{(j,k)\in E^o} x_{ij} = 0 \qquad \forall j \in S, \qquad \text{(A4)}$$

$$T_j \leq T_i + \Delta_i^h + \Delta_i^d + \tau_{ij}^t + M\left(1 - x_{ij}\right) \qquad \forall (i,j) \in E^o \qquad \text{(A5)}$$

$$T_j \geq T_i + \Delta_i^h + \Delta_i^d + \tau_{ij}^t - M\left(1 - x_{ij}\right) \qquad \forall (i,j) \in E^o$$

$$\sum_{k\in H_s}\sum_{(s,j)\in E_s^h} y_{sj}^{sk} - m^h \sum_{(i,s)\in E^o} x_{is} \leq 0 \qquad \forall s \in S, \qquad \text{(A6)}$$

$$\sum_{(s,j)\in E_s^h} y_{sj}^{sk} - \sum_{(i,s)\in E^o} x_{is} \leq 0 \qquad \forall s \in S, k \in H_s, \tag{A7}$$

$$\sum_{(i,j)\in E_s^h} y_{ij}^{sk} - \sum_{(j,l)\in E_s^h} y_{jl}^{sk} = 0 \qquad \forall s \in S, k \in H_s, j \in B, \tag{A8}$$

$$\sum_{s\in S} \sum_{k\in H_s} \sum_{(i,j)\in E_s^h} y_{ij}^{sk} \leq 1 \qquad \forall j \in B, \tag{A9}$$

$$t_j^{sk} \geq t_i^{sk} + \tau_{ij}^h - W^h\left(1 - y_{ij}^{sk}\right) \qquad \forall s \in S, k \in H_s, (i,j) \in E_s^h, \tag{A10}$$

$$t_j^{sk} \leq t_i^{sk} + \tau_{ij}^h + W^h\left(1 - y_{ij}^{sk}\right) \qquad \forall s \in S, k \in H_s, (i,j) \in E_s^h,$$

$$\sum_{(i,j)\in E_s^h} y_{ij}^{sk}\, d_{ij} - \Delta_s^h \leq 0 \qquad \forall s \in S, k \in H_s, \tag{A11}$$

$$\sum_{(s,j)\in E_s^d} z_{sj}^{sk} - \sum_{(i,s)\in E^o} x_{is} \leq 0 \qquad \forall s \in S, k \in D_s, \tag{A12}$$

$$\sum_{k\in D_s} \sum_{(s,j)\in E_s^d} z_{sj}^{sk} - m^d \sum_{(i,s)\in E^o} x_{is} \leq 0 \qquad \forall s \in S, \tag{A13}$$

$$\sum_{(i,j)\in E_s^d} z_{ij}^{sk} - \sum_{(j,l)\in E_s^d} z_{jl}^{sk} = 0 \qquad \forall s \in S, k \in D_s, j \in C, \tag{A14}$$

$$\sum_{s\in S}\sum_{k\in D_s}\sum_{(i,j)\in E_s^d} z_{ij}^{sk} + J_c = 1 \qquad \forall j \in C, \tag{A15}$$

$$t_j^{sk} \geq t_i^{sk} + \tau_{ij}^d - W^d\left(1 - z_{ij}^{sk}\right) \qquad \forall s \in S, k \in D_s, (i,j) \in E_s^d,$$

$$t_j^{sk} \leq t_i^{sk} + \tau_{ij}^d + W^d\left(1 - z_{ij}^{sk}\right) \qquad \forall s \in S, k \in D_s, (i,j) \in E_s^d, \tag{A16}$$

$$\sum_{(i,j)\in E_s^d} z_{ij}^{sk}\, d_{ij} - \Delta_s^d \leq 0 \qquad \forall s \in S, k \in D_s, \tag{A17}$$

$$\sum_{k\in D_s}\sum_{(m,l)\in E_s^d} z_{ml}^{sk} - \sum_{p\in H_s}\sum_{(i,j)\in E_s^h : g_j^l = 1} y_{ij}^{sp} \leq 0 \qquad \forall s \in S, l \in C, \tag{A18}$$

$$\sum_{s\in S}\sum_{k\in D_s} q_j^{sk} + A_c \geq Q_c \qquad \forall j \in C, \tag{A19}$$

$$p_{ij}^{sk} \leq L_{max} z_{ij}^{sk} \qquad \forall s \in S, k \in D_s, (i,j) \in E_s^d, \tag{A20}$$

$$\sum_{j\in C} q_j^{sk} \leq L_{max} \qquad \forall s \in S, k \in D_s, \tag{A21}$$

$$\sum_{(i,j)\in E_s^d} p_{ij}^{sk} - \sum_{(j,l)\in E_s^d} p_{jl}^{sk} = q_j^{sk} \qquad \forall s \in S, k \in D_s, j \in C, \tag{A22}$$

$$\mathrm{T}_{sj} \geq T_s + \Delta_s^h + \sum_{k \in D_s} \mathrm{t}_j^{sk} - M\left(1 - \sum_{k \in D_s} \sum_{(i,j) \in E_s^d} z_{ij}^{sk}\right) \qquad \forall s \in S, j \in C, \tag{A23}$$

$$x_{ij}, y_{kl}^{sf}, z_{mn}^{sp} \in \{0,1\}; J_c, A_c, p_{ij}^{sk}, q_j^{sk}, T_s, \Delta_s^h, \Delta_s^d, t_h^{sf}, \mathrm{t}_c^{sr} \geq 0$$

$$\forall s \in S, (i,j) \in E^o, (k,l) \in E_s^h, (m,n) \in E_s^d, c \in C, h \in B, f \in H_s, r \in D_s, \tag{A24}$$

The objective, shown in (A1), is to minimize a comprehensive cost function, which includes both unfulfillment cost and delay cost components. More specifically, the objective function minimizes: (i) the total delay time to reach all community, and the penalty for (ii) failure to reach communities and (iii) failure to satisfy demand by delivery drones. We omit any fixed cost associated with truck or drone dispatch since we assume that a fixed set of trucks is available, each of them equipped with a fixed number of hotspot and delivery drones.

Constraints (A2)—(A5) restricts the truck routes and satellite visiting times. Constraint (A2) limit the maximum number of trucks that can be used, while the number of trucks allowed to visit a satellite location is restricted to at most one by constraints (A3). Flow preservation for trucks at satellite nodes is ensured by constraints (A4). Constraints (5) dictate the time it takes a truck to reach each satellite node. Constraints (A6)—(A11) dictate the hotspot routes and hotspot location visiting times. In constraints (A6) and (A7), we ensure that hotspot drones are only dispatched from a satellite node if a truck has visited that node; furthermore, the number of hotspot drones dispatched is limited by the number of drones carried by that truck. Constraints (A8) ensure the low preservation for hotspot drones at hotspot location nodes, while the number of hotspot drones allowed to visit a hotspot location is restricted to at most one by constraints (A9). Constraints (A10) dictate the time it takes a hotspot location node by a hotspot drone launched from each satellite. According to constraints (A11), the lengths of time that a truck has to spend at a satellite location due to hotspot drone flights are bounded by the maximum flight times of hotspot routes originating from that satellite, Constraints (A12)—(A17) are analogous to constraints (A6)—(A11), except for they restrict the delivery drone routes and community location visiting times.

Among these constraints, (A15) differs from (A9), which ensure that every community is either visited by a delivery drone or the variable $J_c$ takes a non-zero value. Constraints (A18) are coverage constraints that allow for a community to be visited by a delivery drone dispatched from a satellite only if that location has been covered by a hotspot drone launched from the same satellite. Constraints (A19)—(A22) restrict the load quantities carried and delivered at communities by delivery drones. According to constraints (A19), demand at every community is either satisfied by a delivery drone or the variable $A_c$ takes a non-zero value. Constraints (A20)—(A21) restricts the maximum load quantity that can be delivered by a delivery drone to the drone load carrying capacity, while constraints (A22) ensure the product flow balance at each community node. Constraint set (A23) defines the required times to reach every community by a delivery drone from each satellite. These constraints keep track of the elapsed time between the beginning of the planning period and the time to reach each community by a delivery drone. Finally, variable restrictions are presented in (A24).

### *A2. Comparison between the exact method and CG-based heuristic approach*

We run numerical experiments with small size problem instances to compare the solution quality computational time given by the exact method of solving the MILP model and the CG-bases solution approach we design for the 2EVRP-HD-UAV framework presented in Section 4.

We create datasets for a 20×20 square region, where we randomly select locations for satellites, hotspots, and communities. We create 5 data sets for each of the problem size (with the number of communities 20, 30, 40, and 50). We set the number of satellite locations at 5, number of hotspot locations at 8, and number of trucks at 3. For generating demand and setting other parameter values, we follow the same procedure discussed in Section 6. We summarize our findings in Table A1 below.

Table A1: Comparison between the exact method and column generation-based heuristic

| Problem size | MIP model with exact method | | | | Column generation-based heuristics | | | | Difference between CG-based heuristic and exact method (%) | |
|---|---|---|---|---|---|---|---|---|---|---|
| | CPU time (second) | | Objective function value | | CPU time (second) | | Objective function value | | | |
| | Mean | Standard deviation | Mean | Standard deviation | Mean | Standard deviation | Mean | Standard deviation | Objective function value | CPU time |
| 20 | 1991.0 | 1081.0 | 731.2 | 51.2 | 40.4 | 4.2 | 861.8 | 40.1 | 15.15 | (56.72) |
| 30 | 27794.6 | 2901.4 | 1054 | 27.4 | 107.0 | 39.9 | 1202.8 | 23.4 | 14.12 | (95.67) |
| 40 | * | * | * | * | 152.2 | 16.4 | 1534.2 | 73.7 | - | - |
| 50 | * | * | * | * | 318.0 | 22.3 | 1940.8 | 36.7 | - | - |

* Time limit of 86400 seconds (24 hours) was reached.

### *A3. Pricing subproblems for generating hotspot drone routes (D-PSP-HD-1)*

We define the following two decision variables for these pricing subproblems:

$u_{ij}^{sk}$ =1, if arc $(i,j) \in \bigcup_{s \in S} E_s^h$ is traversed by hotspot drone $k \in H_s$;0 otherwise

$t_j^{sk}$ Time to reach node $j \in B_s \cup s$ by hotspot drone $k \in H_s$

The pricing subproblems resulting from the satellite-drone level decomposition for generating route for each hotspot drone $k \in H_s$ from satellite $s \in S$ are given in (A25)—(A31).

$$Minimize \quad (RC)_{sk}^h = -\gamma_s - \varphi_{sk} + \sum_{(i,j) \in E_s^h} u_{ij}^{sk} \left( \sum_{c \in C_s : g_j^c = 1} \sum_{j \in B_s} \pi_{sc} - \sigma_{sk} \tau_{ij}^h \right) \tag{A25}$$

Subject to:

$$\sum_{(i,j) \in E_s^h} u_{ij}^{sk} \, \tau_{ij}^h \leq W^h \tag{A26}$$

$$\sum_{(s,j)\in E_s^h} u_{sj}^{sk} \leq 1 \tag{A27}$$

$$\sum_{(i,j)\in E_s^h} u_{ij}^{sk} \leq 1 \qquad \forall j \in B_s, \tag{A28}$$

$$\sum_{(i,j)\in E_s^h} u_{ij}^{sk} - \sum_{(j,i)\in E_s^h} u_{ij}^{sk} = 0 \qquad \forall j \in B_s \tag{A29}$$

$$\begin{aligned} t_j^{sk} &\geq t_i^{sk} + \tau_{ij}^h - W^h\left(1 - u_{ij}^{sk}\right) && \forall (i,j) \in E_s^h, \\ t_j^{sk} &\leq t_i^{sk} + \tau_{ij}^h + W^h\left(1 - u_{ij}^{sk}\right) && \forall (i,j) \in E_s^h, \end{aligned} \tag{A30}$$

$$u_{ij}^{sk} \in \{0,1\}, \qquad (i,j) \in E_s^h; \qquad t_j^{sk} \geq 0, \qquad \forall j \in B_s, \tag{A31}$$

The objective function in (A25) finds the best reduced cost $(RC)_{sk}^h$ of a route for hotspot drone $k \in H_s$ launched from satellite $s \in S$. If the obtained reduced cost is negative, then that route is added to the restricted set $U_s$. Constraints (A26) limits the total travel time of a drone. Each drone can be used only along a single path starting from the satellite under consideration according to constraints (A27). Each hotspot location can be visited at most once per constraints (A28); constraints (A29) provide the balancing of inbound and outbound arcs coincident to a hotspot location. Constraints (A30) indicate the time to visit two consecutive nodes on hotspot drone route. Finally, constraints (A31) provide non-negativity bounds. As a reminder, $\gamma_s, \varphi_{sk}, \sigma_{sk}, \pi_{sc}; \forall s \in S, k \in H_s, j \in C_s$ in (A25) are the dual variables associated with constraints (6), (8), (11) and (14) in the D-RMP-1 model presented in subsection 4.1.

### *A4. Pricing subproblems for generating routes in satellite-level decomposition*

The pricing subproblems for generating routes for hotspot drones launched from satellite $s \in S$, are given in (A32)—(A33).

$$Minimize \quad \sum_{k \in H_s} \left( -\gamma_s - \varphi_{sk} + \sum_{(i,j) \in E_s^h} u_{ij}^{sk} \left( \sum_{c \in C_s : g_j^c = 1} \sum_{j \in B_s} \pi_{sc} - \sigma_{sk} \tau_{ij}^h \right) \right) \tag{A32}$$

Subject to:

$$(A26) - A(31) \tag{A33}$$

In this pricing subproblems, we impose the constraints (A26)—(A31) for each hotspot drone $k \in H_s$ from satellite $s \in S$.

### *A5. Pricing subproblems for generating delivery routes in satellite-level decomposition*

The pricing subproblems for generating route for delivery drones launched from satellite $s \in S$, are given in (A33)—(A34).

$$\begin{aligned} Minimize \quad \sum_{k \in D_s} \Bigg( -\theta_s - \phi_{sk} - \sum_{j \in C_s} \sum_{(i,j) \in E_s^d} v_{ij}^{sk} \left( \zeta_j + \pi_{sj} + \rho_{sk} \tau_{ij}^d \right) \\ - \sum_{j \in C_s} \alpha_{sj} \left( t_j^{sk} + M \sum_{(i,j) \in E_s^d} v_{ij}^{sk} \right) \Bigg) \end{aligned} \tag{A33}$$

Subject to:

$$(17) - (22) \tag{A34}$$

In this pricing subproblems, we impose the constraints (17)—(22) for each delivery drone $k \in D_s$ launching from satellite $s \in S$.